# Robot Path Planning by Traveling Salesman Problem with Circle Neighborhood: modeling, algorithm, and applications

Arman Nedjati[a], Béla Vizvári[b]

[a]Industrial Engineering Department, Quchan University of Technology, Quchan, Iran

[b]Industrial Engineering Department, Eastern Mediterranean University, Famagusta, Via Mersin10, North Cyprus

**Abstract**

The study investigates the problem of traveling salesman problem with circular neighborhood (TSPCN). Instead of cities there are circles and each point on circle can be a potential visiting node. The problem is to find the minimum length Hamiltonian cycle connecting the circles. Among the various real life applications of the problem, this paper concentrates on robot path planning for the laser welding robot, and data collection in a wireless sensor network by unmanned aerial vehicles (UAVs). The TSPCN is formulated as a nonlinear model, the objective function is linearized, and as a solution procedure decomposed into a two-phase model. The Models are coded in Cplex and Knitro and solved for small and medium sized instances.

## 1. Introduction

Traveling Salesman Problem (TSP) is a well-known problem in the field of optimization. Plenty of variants of TSP among with mathematical models and heuristic algorithms are applied in real life situations[1–3]. TSP has many applications in transportation, robotics, manufacturing, computer networks, etc. One of the important applications that has many generalizations is path planning in robotics [4–8]. This study focuses on the problem of path planning for a mobile robot (i.e. unmanned aerial vehicle (UAV) or a welding robot arm) that has to traverse the minimum path to connect some predefined circles.

A real life problem that was a motivation for this study is related to Wireless Sensor Networks (WSN) that have been widely investigated and improved in last decade. The storage, communication, processing, and sensing capabilities of sensors are increased for agricultural, boarder monitoring, military uses, etc. When the sensors are distributed in a widespread area that is not possible to establish a communication network, a UAV could traverse to collect the data by passing through the spherical radio range around the sensors. However UAV path planning has many applications [9], WSN interested readers can refer to earlier studies [10–14].

Moreover another industrialized application of proposed problem is related to welding robot path planning. The spot welding robots are extensively used in processing industry, machinery, aerospace, and automotive industries. These robots are effective when the number of welding points are increased. An optimal path for welding robot arm can effectively shorten the overall process time and consequently improves the productivity[15]. Both assembling and laser welding are typical examples. Industrial robots are still expensive equipment. A laser ray has high power, therefore the walls of the cell where the robot works must be strengthened by thick steel sheets to obtain a secure environment. It also increases the costs further. In this application a new optimization method is suggested for designing the route and to solve this problem. For another approach of designing the route of the laser robot arm see [1].

According to numerous real life applications of TSP, many variants [1,16] of this problem such as generalized TSP (GTSP) [17], TSP with clustering [18], hierarchical TSP [19], colored TSP [20], etc. have been investigated through previous decades. Moreover, there are many studies that focused on robotics applications of TSP. [21] addresses the problem of mobile onsite power delivery for autonomous underwater vehicles as a multiple GTSP to improve the range and duration of underwater operations. As a collaborative human-robot TSP [22] developed a model based on an initial remote sensing to optimize the selective TSP tour design and seeking the human assistance. Even forest management activities by robots are also investigated under the concept of TSP [23]. As a comprehensive survey [24] concentrated on a broad review of the classical studies and latest breakthroughs concerning the TSP and routing approaches applied in robotics. Furthermore a wide range of studies investigate on different types of TSP algorithms consists of approximation [25], heuristic [26], meta-heuristic [27], hybrid [28], and exact solution procedures [29].

In this study we present a new industrialized variant of TSP which is TSP with Circle Neighborhood (TSPCN). TSPCN is a variant of Traveling Salesman Problem with Neighborhoods (TSPN). TSPN was first introduced by [30] and investigated in approximation algorithms [31]. In TSPN n number of connected neighborhoods or regions (as subsets) from an Euclidean plane, and a set of predefined coordinates are given, and the problem is to minimize the length of the tour that visits each neighborhood in one of the predefined coordinates [31]. Another relevant problem is GTSP where there are groups of coordinates within number of regions and the task is to find the minimum length Hamiltonian path that visits any coordinate of each group exactly once. In fact the TSPN is a generalization of GTSP. Therefore, since the TSPCN problem generalizes the Euclidean TSP when the circle radius is zero, it is considered as a NP-hard problem. As in TSPCN the regions are the points located on any coordinate of a circle and the tour can visit each region at any of its coordinate the problem is categorized as continuous TSPN [31]. In TSPCN there are N number of circles (cities) that each must be visited once while visiting node can be any point inside or on the perimeter of the circles. Regarding the literature the TSPCN problem has not been studied before in this form. Fig 1 depicts a feasible solution of problem for a 10 circuit problem.

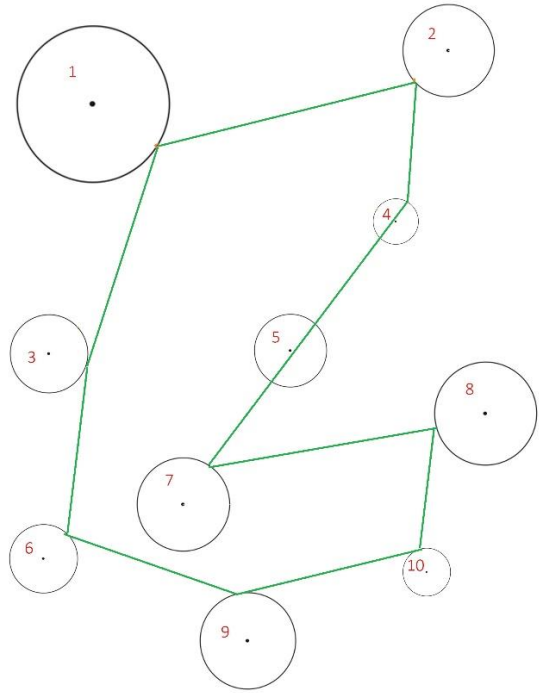

Figure 1. A feasible solution of Ci-TSP

The remainder of the article is structured as follows: The nonlinear mathematical formulations among with some linearization and a two phase approach are presented in Section 2. Some small and medium scale instances are tested and solved by Baron, Couenne, and Knitro nonlinear solvers in Section 3. Discussion, conclusion and some suggestions for future studies are offered in section 4.

## 2. Mathematical Formulations

### 2.1 Nonlinear mathematical model

A horizontal plane contains $N$ points. Each point is an extreme point of a cone. The intersections of the cones with another horizontal plane above the first one, are the circles $C_1, C_2, \ldots, C_n$. In the case of the laser welding robot, the surface of the product has a safety cover with the exception of the welding points. The cone of a point is the set of points of the space from where the laser ray can "see" the welding points. Similarly, the cone of a wireless sensor consists of the points of the space where the broadcast of the sensor can be detected. There is no restriction on the radius or position of the circles. The problem is to select one point from each circle and construct a tour along the selected points with a minimal total length. The selected points of different circles may coincide. In the latter case the distances of the cities having the same selected point is zero.

Let $C_i = \{y_i = (u_i, v_i) | (u_i - a_i)^2 + (v_i - b_i)^2 \leq r_i^2\}$, where $(a_i, b_i)$ is the center of the circle and $r_i$ is the radius.

The well-known Miller, Tucker and Zemlin formulation of TSP [7] is applied in our model. The notations used in the model below are as follows:

**Sets**:

C :{1,2,…,N} the index set of circuits

**Parameters**:

$N$ = the number of circles.
$(a_i, b_i)$ = the center of circle $i$.
$R_i$ = the radius of circle $i$.

**Variables**:

$seq_i$ : Sequence in which circle $i$ ($i \neq 1$) is visited.
$(u_i, v_i)$ : The selected point of circle $i$
$dis_{ij}$ : The distance of the selected points of circles $i$ and $j$
$p_{ij}$ : Binary variable; equals to one if circle $j$ is visited immediately after circle $i$, otherwise it is zero.

The objective function is to minimize the total distance.

Model $T_0$:

$$Min \sum_{i \in C} \sum_{j \in C} p_{ij} dis_{ij} \tag{1}$$

St;

$$\sum_{j \in C} p_{ij} = 1 \quad \forall i \in C | i \neq j \tag{2}$$

$$\sum_{j \in C} p_{ji} = 1 \quad \forall i \in C | i \neq j \tag{3}$$

Constraints (2) and (3) ensures that each node must be visited exactly once.

$$seq_i - seq_j + Np_{ij} \leq N - 1 \quad \forall i \neq j \in C | i, j \neq 1 \tag{4}$$

Constraint (4) is subtour elimination constraint adopted from [32].

$$R_i^2 - (u_i - a_i)^2 - (v_i - b_i)^2 \geq 0 \quad \forall i \in C \tag{5}$$

Constraint (5) guarantees that the selected nodes must be in the circuit.

$$dis_{ij}^2 \geq (u_i - u_j)^2 + (v_i - v_j)^2 \quad \forall i \neq j \in C \tag{6}$$

$$p_{ij} \in \{0,1\} \quad \forall i,j \in C \tag{7}$$

Constraint (6) assures that the value of variable $dis_{ij}$ cannot be less than the distance of the two selected points, and constraint (7) represents the integrality constraint. The optimality will ensure strict equation if it is necessary.

As the objective function is non-linear in its original form (1), it is possible to linearize it. It is necessary to rewrite the objective function and constraint (6) as follows:

$$Min \sum_{i \in C} \sum_{j \in C} dis_{ij} \tag{8}$$

$$dis_{ij}^2 + M(1 - p_{ij}) \geq (u_i - u_j)^2 + (v_i - v_j)^2 \quad \forall i \neq j \in C \tag{6-a}$$

$$dis_{ij} \geq 0 \tag{6-b}$$

Where M is a very big number.

### 2.2 Two phase mathematical formulations

It is important to emphasize in connection to the numerical solution that it is an industrial design problem as the route of the robot arm is optimized only once but the same solution is applied many times. The number of applications is at least the number of cars produced from the car model. Thus, it is worth to invest a longer CPU time into the solution of the mathematical model to obtain a good feasible solution for the industrial case. However, in case of data collection application (WSN and UAVs), CPU time is important. Thus, in order to fasten the solution procedure a two-phase model is presented in the following. In the first phase, the circles are divided into four equal quarters. Each quarter is substituted by the middle point of its arc. Only one out of these four points needs to be visited. The optimal solution from the first phase determines a sector for each circle and any of the points of this sector can be selected in the second phase.

**Sets**:

L :{1,2,3,4} Set of nodes on circuits

**Parameters:**

$H_{ij}^{kl}$ = Distance between node $k$ located on circuit $i$ and node $l$ located on circuit $j$.

$G_l$ = The angle formed by the radius, with respect to the positive X-axis (degrees).

**Variables:**

$x_{ij}^{kl}$ : Binary variable equals to 1 if there is an arc from node $k$ of circuit $i$ to node $l$ of circuit $j$, otherwise it is 0.

$$H_{ij}^{kl} = \sqrt{((R_i Cos(G_l) + a_i) - (R_j Cos(G_k) + a_j))^2 + ((R_i Sin(G_l) + b_i) - (R_j Sin(G_k) + b_j))^2} \quad (11)$$

Equation (11) which is depicted in Fig 2 calculates the distance between node $l$ of circuit $i$ and node $k$ of circuit $j$.

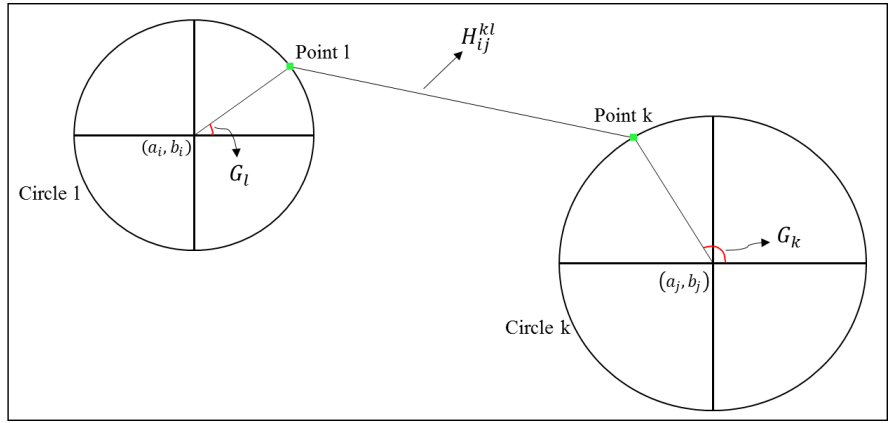

Figure 2. The explanation of Formula (11)

Model $T_1$ is the first phase that finds the optimal tour while each circuit's perimeter is divided into quarters and the center of each quarter is a potential node. Thus model $T_1$ finds the optimal sequence of visiting the circuits based on the center point of each quarter's perimeter. Then Model $T_2$ uses the optimal visiting sequence as a parameter to find the optimal coordinates.

Model $T_1$:

$$Min \sum_{i \in C} \sum_{j \neq i \in C} \sum_{l \ni L} \sum_{k \in L} H_{ij}^{kl} x_{ij}^{kl} \quad (12)$$

St:

$$\sum_{j \in C} \sum_{l \in L} \sum_{k \in L} x_{ij}^{kl} = 1 \quad \forall i \in C \quad (13)$$

$$\sum_{j \in C} \sum_{l \in L} \sum_{k \in L} x_{ji}^{kl} = 1 \quad \forall i \in C \quad (14)$$

$$\sum_{l \in L} \sum_{k \in K} x_{ij}^{kl} + \sum_{l \in L} \sum_{k \in K} x_{ji}^{lk} \leq 1 \quad \forall i \neq j \in C \quad (15)$$

$$seq_i - seq_j + N \sum_{l \in L} \sum_{k \in L} x_{ij}^{kl} \leq N - 1 \quad \forall i \neq j \in C \,|\, i,j \neq 1 \quad (16)$$

$$\sum_{j \in C} \sum_{k \in L} x_{ij}^{lk} - \sum_{j \in C} \sum_{k \in L} x_{ji}^{kl} = 0 \quad \forall i \in C, l \in L \tag{17}$$

$$x_{ij}^{kl} \in \{0,1\} \quad \forall i,j \in C, \forall l,k \in L \tag{18}$$

$$seq_i \in integer \ and \ \geq 1 \quad \forall i \in C \tag{19}$$

Objective function (12) minimizes the total distance of the tour and the constraint are well known TSP problem equations. However constraint (15) is a valid inequality and helps to remove undesirable fractional solutions and tighten the relaxation of the model.

Let's define additional sets of parameters for model T2 as follows;

$O_{il}$: Equals to one if node $i$ of circuit $l$ is visited, and zero otherwise.

$q_{ij}$: Equals to one if we move to circuit $j$ right after circuit $i$.

$Umax_i$: Maximum amount on X-axis that the optimal point can take on circuit $i$.

$Umin_i$: Minimum amount on X-axis that the optimal point can take on circuit $i$.

$Vmax_i$: Maximum amount on Y-axis that the optimal point can take on circuit $i$.

$Vmax_i$: Minimum amount on X-axis that the optimal point can take on circuit $i$.

$$q_{ij} = \sum_{l \in L} \sum_{k \in L} x_{ij}^{kl} \quad \forall i,j \in C \tag{20}$$

Equations (20) and (21) assigns the results of model T1 to the parameters that will use in model $T_2$.

Model $T_2$:

$$Min \sum_{i \in C} \sum_{j \in C} q_{ij} dis_{ij} \tag{21}$$

St:
$$eq \ (5), (6)$$
$$Umin_i \leq u_i \leq Umax_i \quad \forall i \in C \tag{22}$$

$$Vmin_i \leq v_i \leq Vmax_i \quad \forall i \in C \tag{23}$$

Objective function (21) minimizes the total distance while constraints (22) and (23) assigns the lower bound and upper bound to coordinate variables based on the results of model $T_1$.

## 3. Results

Two examples for problems with 12 and 20 circles can be seen in Figure 3 and Figure 4. The optimal solution was obtained by the two phase method described in section 2.2. The main nonlinear model was coded in GAMS and solved by different nonlinear solvers (Couenne, Antigone, Baron, Knitro), however no integer solution was obtained. The two-phase method solved the 12-circle's problem optimally in 500 seconds and for 20-circle problem in 2400 seconds. The optimal solution is shown in Figure 3.

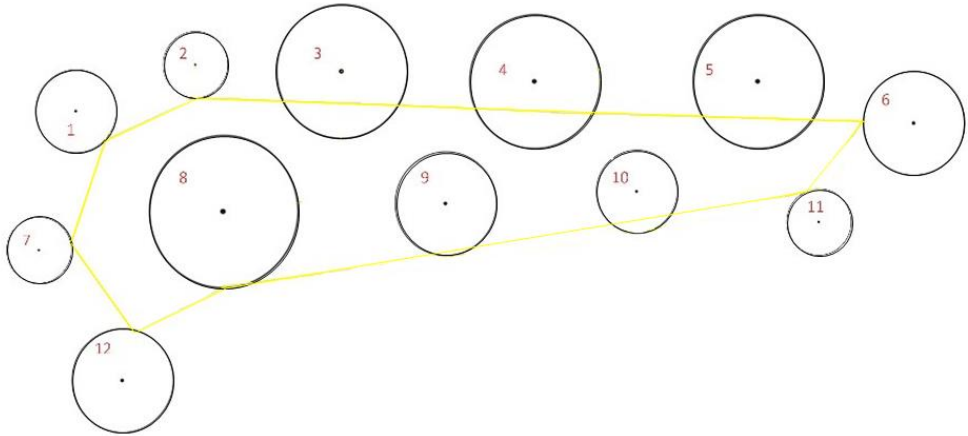

Figure 3 Optimal solution of 12 circuits (Two-phase method)

Figure 4. Optimal solution of a 20-circle problem (Two-phase method)

Also a 40 circuit problem is solved by the two phase method with the execution time of 6700 seconds. Figure 4 depicts the optimal solution of the problem.

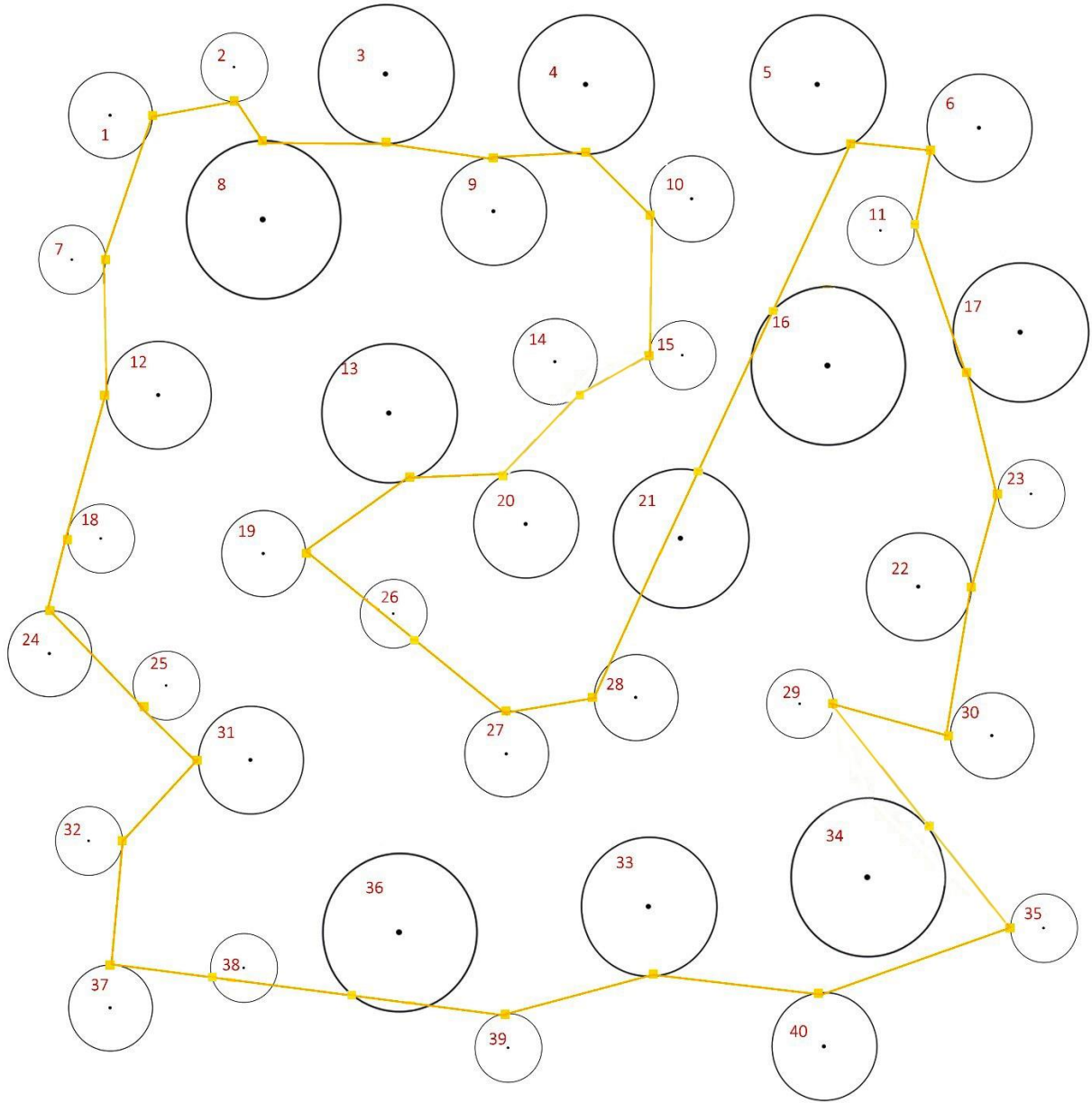

Figure 5 The optimal solution of a 40 circuit problem (Two-phase method)

Another problem is solved. This time it has a real industrial size of 75 welding points. The problem and its optimal solution can be found on Figure 5. The original mathematical model was too difficult in this case. Therefore the problem was solved in two steps. In the first step a traditional TSP was solved. The positions of the cities were the center points of the regions. The selected points $(u_i, v_i)$ were also included in the problem for some technical reasons although they did not affect the optimality in any sense. In the second main step those subtour elimination constraints were used which were generated in the first main step.

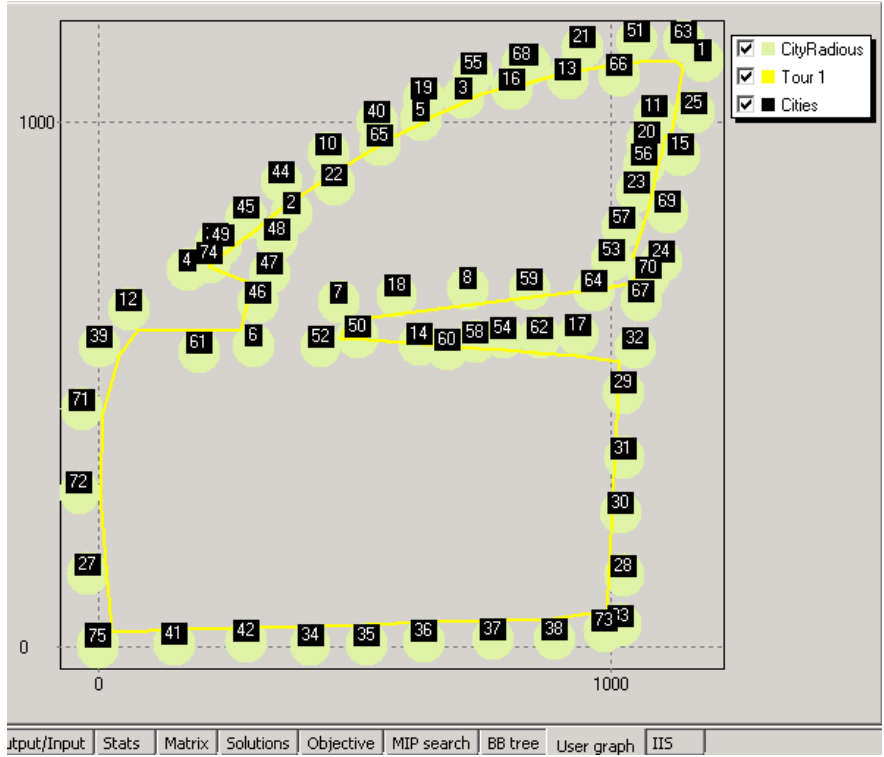

Figure 6 Industrialized size problem of 75 circuits

The method mentioned above is a so-called cutting plane method. As the problem is too large to solve, instead of the original problem a relaxation is solved. In many cases, relaxation means that some constraints are deleted. It means that not all subtour elimination constraints are added to the problem. If the optimal solution of the relaxed problem does not satisfy all the constraints then a violated subtour elimination constraint is added to the problem. The name "cutting plane method" reflects the fact that the subtour elimination constraint is a linear inequality which cuts the previous optimal solution.

## 4. Discussion and conclusion

This short discussion paper investigates a new problem of robotic path planning. The study concentrates on UAVs path that collect data from distributed wireless sensors, and welding robot path in automobile industry. The problem is a generalization of traveling salesman problem with neighborhood when the neighborhood shapes are circles with different radiuses. It is enough to visit just one point of the circle while the minimum length Hamiltonian cycle. This problem could also have applications in telecommunication systems. The exact nonlinear mathematical model of the problem was developed and a two phase method along with a linearization are presented to shorten the execution time. The models are validated with some numerical examples and the results shows the efficiency of the two-phase method.

Direction for future studies can be categorized in the subsequent three fields: Heuristic and Meta-heuristics development, problem and modeling extension, and exact solution procedures. Multi-robot system or visiting time window for sensors could be some examples of problem extension. Moreover, instead of a circle, one can consider the neighbors as a cut of the hemisphere, while the complexity of problem is arisen, the accuracy is increased.